\documentclass[12pt]{article}
\usepackage{amssymb, cite,graphicx}

\def\ben{\begin{enumerate}}
\def\een{\end{enumerate}}

\newtheorem{stat}{}[section]
\def\bs{\begin{stat}}
\def\es{\end{stat}}

\def\bp{\noindent{\bf Proof}  \ \ \ }
\newcommand{\ep}{\hfill $\square$}

\makeatletter
\@addtoreset{equation}{section}

\makeatother

\begin{document}
\begin{center}
{\Large {{\bf PACKING $k$-EDGE TREES IN GRAPHS}}}
\\[2ex]
{\Large {{\bf OF RESTRICTED VERTEX DEGREES}}}
\\[4ex]
{\large {\bf Alexander K. Kelmans}}
\\[2ex]
{\bf University of Puerto Rico, San Juan, Puerto Rico}
\\[0.5ex]
{\bf Rutgers University, New Brunswick, New Jersey}
\end{center}

\begin{abstract}      
Let ${\cal G}^s_r $ denote the set of graphs with each vertex of degree at least
$r$ and at most $s$,  $v(G)$ the number of vertices, and 
$\tau _k (G)$ the maximum number of disjoint $k$-edge trees in $G$.
In this paper we show that 
\\[0.5ex]
$(a1)$ if $G \in {\cal G}^s_2 $ and $s \ge 4$, then $\tau _2 (G) \ge v(G)/(s+1)$,
\\[0.5ex]
$(a2)$ if $G \in {\cal G}^3_2 $ and $G$ has no 5-vertex components, then
$\tau _2 (G) \ge v(G)/4$,
\\[0.5ex]
$(a3)$ if $G \in {\cal G}^s_1 $ and $G$ has no $k$-vertex component, where $k \ge 2$ and $s \ge 3$, 
then
$\tau _k(G) \ge (v(G) - k)/(sk - k +1)$, and
\\[0.5ex]
$(a4)$ the above bounds are attained for infinitely many connected graphs.
\\[0.5ex]
\indent
Our proofs provide polynomial time algorithms for finding 
the corresponding packings in a graph. 
\\  

\end{abstract}

\section{Introduction}
\label{intro}

\indent

We consider simple undirected graphs.  
All notions on graphs that are not defined here can be found in  \cite{BM,D}.  
Given a graph $G$ and a set ${\cal F}$ of subgraphs of $G$, 
an ${\cal F}$-{\em packing} of $G$ is a subgraph of $G$ whose components are 
members of ${\cal F}$. 
The ${\cal F}$-{\em packing problem} is that of
finding an ${\cal F}$-packing having the maximum number of vertices. 
Various ${\cal F}$-packing problems
have extensively  been studied by many authors for different families
${\cal F}$  (e.g. \cite{CH,EKK,HK,KKN,K,KMS,KM,LV,LoPo,LoPo1,LP}).

It is not surprising that ${\cal F}$-packing problem turns out 
to be $NP$-hard for most of the families ${\cal F}$.
Surprisingly  the problem can be solved in polynomial time for some 
non-trivial families. 
For example, Edmonds (see \cite{LP}) showed that
a classical matching problem can be solved in polynomial time. 
It is also known that the problem of packing stars of at least one
and at most $k$ edges is polynomially solvable even if the stars, we pack,
are required to be induced \cite{K,K97}.
On the other hand, the problem of finding in $G$ the maximum number of
disjoint subgraphs, isomorphic to a given connected graph $H$ of at least
three vertices, is $NP$-hard \cite{HK}. 
The problem remains $NP$-hard for cubic graphs if $H$ is 
a path having at least two edges \cite {Khrd}.

In this paper we consider the $H$-packing problem when $H$ is a tree and,
in particular,  when $H = \Lambda $, a path of two edges. 
Although the $\Lambda $-packing problem is $NP$-hard, i.e.
possibly intractable in general, this problem turns out to be tractable for
some natural classes of graphs. 
Here are some examples of such results.

Let $v(G)$ denote the number of vertices of a graph $G$ and
$\tau _k (G)$ denote the maximum number of disjoint 
$k$-edge trees in $G$.
We also put $\tau _2 (G) = \lambda (G)$. 

A graph is called {\em claw-free} if it contains no induced subgraph 
isomorphic to $K_{1,3}$ (which is called a {\em claw}). 

A block $B$ of a  graph $G$ is called an {\em end-block} of $G$ if $B$ has  exactly 
one vertex adjacent to a vertex in $G - B$.

Obviously $\lambda (G) \le \lfloor v(G)/3 \rfloor $.
\bs {\em \cite{KKN}}
\label{2conclfr} 
Suppose that $G$ is a connected claw-free graph having at most two end-blocks 
$($in particular, a 2-connected claw-free graph$)$.
Then $\lambda(G) = \lfloor v(G)/3 \rfloor $.
\es 

Let $G^\Delta $ be the graph obtained from a cubic graph by 
``replacing'' every vertex by a triangle. 
From {\bf \ref{2conclfr}} we have, in particular, for cubic graphs:

\bs {\em \cite{KKN}}
\label{2contrcb} 
Suppose that $G$ is a connected cubic graph having at most two end-blocks 
$($in particular, a 2-connected graph$)$.  
Then
$\lambda(G^\Delta ) = \lfloor v(G^\Delta )/3 \rfloor $.
\es 

Let $eb(G)$ denote the number of end-blocks of $G$.
The previous theorem follows from the following more general result.
\bs {\em \cite{KKN}}
\label{eb(G)clfr}
Suppose that $G$ is a simple connected claw-free graph and  
$eb(G) \ge 2$. 
Then  
$\lambda(G) \ge \lfloor (v(G) - eb(G) + 2)/3\rfloor$,
and this lower bound is sharp.
\es

\bs {\em \cite{KMS}}
\label{treepack}
Let $T$ be a tree on $t$ vertices and
let $\epsilon > 0$.
Suppose that $G$ is a $d$-regular graph on $n$ vertices
and $d \geq \frac{128t^3}{\epsilon^2}\ln( \frac{128t^3}{\epsilon^2})$.
Then $G$ contains at least $(1-\epsilon)n/t$ vertex disjoint copies of $T$.
In particular, if $G_n$ is a $d_n$-regular graph on $n$ vertices
and $d_n \to \infty$ when $n \to \infty$, 
then $G_n$ contains at least $(1-o(1))n/t$ 
$($and, obviously, at most $n/t$$)$ disjoint trees isomorphic to $T$.
\es

\bs  {\em \cite{KM}}
\label{1}  
Let $G$ be a cubic graph. Then $\lambda (G) \ge v(G)/4$.  
\es 

Let ${\cal G}^s_r $ denote the set of graphs with each vertex of degree at least
$r$ and at most $s$.
Our main question on $\Lambda $-packings is:
\\[1ex]
\indent
{\em How many disjoint 2-edge paths must an $n$-vertex graph $G$ from 
${\cal G}^s_2$ have?}
\\[1ex]
\indent
In  {\bf \ref{Main1}} and {\bf \ref{DMain1}} we give corresponding lower bounds
on the numbers in question. 
We also show  (see {\bf \ref{extrgraphs1}})
that these bounds are tight.
\bs
\label{Main1} Suppose that $G \in {\cal G}^3_2$ and  
$G$ has no 5-vertex components.
Then $\lambda (G) \ge v(G)/4$.
\es

Obviously {\bf \ref{1}}  follows from {\bf \ref{Main1}} because
if $G$ is a cubic graph then $G \in {\cal G}^3_2$ and $G$ has no 5-vertex
components. 
\\

\bs
\label{DMain1}
 Let $G \in {\cal G}^s_2$, $s \ge 4$.
Then $\lambda (G) \ge v(G)/(s+1)$.
\es

We  give constructions (see Section \ref{constr}) that allow us to prove the following:
\bs
\label{extrgraphs1}
There are infinitely many connected graphs for  which the bounds  
in {\bf \ref{Main1}} and {\bf \ref{DMain1}} are attained.
Moreover, for every integer $k$, such that $1 \le k \le s$, 
there are infinitely many subdivisions of
$k$-regular graphs for which the bounds in 
{\bf \ref{Main1}} and {\bf \ref{DMain1}} are attained.
\es

We also consider a special case of the ${\cal F}$-packing
problem when ${\cal F}$ is the set of all connected subgraphs of $G$ having $k$ edges. 
Let ${\cal T}_k$ denote the set of all trees having $k$ edges.
Obviously this problem is equivalent to the 
${\cal T}_k$-packing problem, the problem of finding in a graph 
the maximum number of disjoint trees of $k$ edges.
Notice that Theorem {\bf \ref{treepack}} provides a similar asymptotic 
result for the ${\cal T}_k$-packing problem.
Our question on ${\cal T}_k$-packings is:
\\[1ex]
\indent
{\em How many disjoint $k$-edge trees must an $n$-vertex graph $G$ from 
${\cal G}^s_1$ have?}
\\[1ex]
\indent
In  {\bf \ref{main1}} we give a lower bound on the number in question. 
We also show  (see {\bf \ref{extrgraphs2}}) that these bounds are tight.

Our main result on the ${\cal T}_k$-packing problem are the following:
\bs
\label{main1} Let $s$ and $k$ be integers, $s \ge 3$, $k \ge 2$.
Suppose that $G \in {\cal G}^s_1$ and $G$ has no $k$-vertex component. 
Then $\tau _k(G) \ge (v(G) - k)/(sk - k + 1)$.
\es

One of the constructions in Section \ref{constr} allows to prove:
\bs
\label{extrgraphs2}
There are infinitely many connected graphs for  which the bound  
in {\bf \ref{main1}} is attained.
\es

Our proofs provide polynomial time algorithms for finding 
the corresponding packings in graphs from ${\cal G}_r^s$. 
Thus these algorithms are polynomial approximation algorithms 
for the corresponding problems.
\\

The results of this paper were presented at the Workshop ``Graph Partitions'' 
in DIMACS, Rutgers University, in July, 2000
(see also  \cite{Kdmx}).

\section{Main notions, notation and simple observations}

\indent

Let $G$ be a simple graph.    
We use the following notation and notions:  
\\[1ex]  
$V = V(G)$ and $E = E(G)$ are the sets of vertices and edges of a graph $G$,  
respectively, 
\\[0.7ex]  
$v(G) = |V(G)|$ and $e(G) = |E(G)|$,  
\\[0.7ex] 
if $x \in V(G)$, then $d(x,G) = d(x)$ is the degree of $x$ in $G$,
\\[0.7ex]
$xPy$ is a path with the end-vertices $x$ and $y$,
\\[0.7ex]
an ${\cal F}$-{\em packing} of $G$ is a subgraph of $G$ whose components are 
members of ${\cal F}$,
\\[0.7ex]
a $\Lambda $-{\em packing} is a subgraph of $G$ whose components are 
2-edge paths in $G$,
\\[0.7ex]
$\lambda (G)$ is the maximum number of vertex disjoint
2-paths in $G$,
\\[0.7ex]
$\tau _k (G)$ is the maximum number of disjoint $k$-edge trees in $G$,
where $k \ge 1$, and so $\lambda (G) = \tau _2(G)$.
\\[1ex]
\indent
A {\em path-thread} in a graph $G$ is a maximal path $P$ in $G$ such that
each vertex of degree two in $P$ is also of degree two in $G$.
A {\em cycle-thread} in $G$ is a cycle $C$ in $G$ such that each vertex of $C$,
except for 
one,  is of degree two in $G$.
A {\em thread} is either a path-thread or a cycle-thread.
\\[1ex]
\indent
A {\em block} of a connected graph $G$ is a maximal connected subgraph $H$ of $G$ such that $H - v$ is connected for every $v \in V(H)$. 
A block $B$ of a 
graph $G$ is called an {\em end-block} of $G$ if $B$ has  exactly 
one vertex adjacent to a vertex in $G - B$.
\\[1ex]
\indent
A  {\em leaf} of a graph $G$ is either a vertex of degree one or an end-block of  at least two edges in $G$.
An $H$-{\em leaf} is a leaf isomorphic to a graph $H$.
If $k$ is an integer, then a $k$-{\em leaf} is a leaf having $k$ vertices. 
\\[1ex]
\indent
Let, as above, ${\cal G}^s_r$ denote the set of graphs with 
each vertex of degree at least $r$ and at most $s$.
\\[1ex]
\indent
Let ${\cal F}^3_2$ denote the set of graphs such that
\\[1ex]
$(a1)$ $d(x,G) \in \{2,3\}$ for every vertex $x$ of $G$ and
\\[1ex]
$(a2)$ $G$ has no 5-vertex components.
\\[1ex]
\indent
Given a class of graphs ${\cal A}$, a graph $G$ is called 
${\cal A}$-{\em mininimal} if $G \in {\cal A}$ and 
$G - e \not \in {\cal A}$ for every $e \in E(G)$. 
Obviously ${\cal G}^s_r$-minimal and 
${\cal F}^3_2$-minimal graphs exist.
\\

It is easy to see the following.
\bs
\label{min}
Let $G$ be a  graph.
\\[1ex]
$(c1)$ If $G \in {\cal F}^3_2$ $($$G \in {\cal G}^s_r$$)$, then $G$ has an ${\cal F}^3_2$-minimal 
$($respectively,  ${\cal G}^s_r$-minimal$)$  
spanning subgraph $F$ and $\lambda (G) \ge \lambda (F)$.
\\[1ex]
$(c2)$ A graph $G$ is ${\cal F}^3_2$-minimal if and only if each edge of $G$ is incident to either a vertex of degree two in $G$ or to a vertex of a $5$-leaf of $G$.
\es

\bs
\label{leaf}
Let $G \in {\cal F}^3_2$, $Q$ be a leaf of $G$.
Then $Q \in {\cal F}^3_2$, the boundary vertex $q$ of $Q$ is of degree 
two in $Q$, and there is a $($unique$)$ path-thread $xTy = S(Q)$ in $G$ 
such that $x = q \in V(Q)$ and $y \in V(G - Q)$.
\es

If $Q$ is a leaf of a graph $G \in {\cal F}^3_2$, 
then $S(Q)$ is called the {\em stem} of the leaf $Q$.
\\ 
Let $\grave{Q}$ be the graph obtained from $Q \cup S(Q)$
by removing the vertex of degree one.

\section{Graph reductions}
\label{reductions}
\indent

Our proving strategy is to establish various properties of 
a minimum counterexample, and finally conclude that it cannot exist. 
At various stages in this process, we will need some operations that reduce 
a minimum counterexample to a smaller one providing a contradiction.  
In this section we present such reductions.
\bs
\label{Dthread1} 
Let $G$ be a graph, $xTy$ be a thread of $G$, where possibly $x = y$,
$T' = T -\{x,y\}$, and $G' = G - T'$
{\em (see Figure \ref{fig:1a})}.
Let $s \ge 3$ be an integer.
\\[1ex]
$(a1)$ 
If $s \ge 4$ and $e(T) \ge 4$, then 
$\lambda (G') \ge v(G')/(s+1) \Rightarrow \lambda (G) \ge v(G)/(s+1)$.
\\[1ex]
$(a2)$
If $s = 3$ and $v(T') \in \{3k, 3k+1\}$, where $k \ge 1$ is an integer, then
\\
$\lambda (G') \ge v(G')/4 \Rightarrow \lambda (G) \ge v(G)/4$.
\es
\bp 
Clearly $v(G) = v(G') + v(T')$ and 
$\lambda (G) \ge \lambda (G') + \lfloor v(T')/3 \rfloor$.
Now $\lambda (G') \ge v(G')/(s+1) 
\Rightarrow 
\lambda (G) \ge v(G')/(s+1) + \lfloor v(T')/3 \rfloor =
v(G)/(s+1) + \gamma _s(T')$
where $\gamma _s(T') = \lfloor v(T')/3 \rfloor - v(T')/(s+1)$.
Since $e(T) \ge 4$, clearly $v(T') \ge 3$.

Suppose that $s \ge 4$. Then obviously $\gamma _s(T') \ge 0$.

Now suppose that $s = 3$. 
If $v(T') = 3k$, then $\gamma _3(T') = k - 3k/4 = t/4$.
If $v(T') = 3k+1$, 
then $\gamma _3(T') = k - (3k+1)/4 = (k-1)/4 \ge 0$.
Since $k \ge 1$,  in both cases 
$\lambda (G) \ge v(G)/(s+1)$.
\ep

\bs
\label{thread2} 
Let $G$ be a graph, $xTy$ be a thread of $G$,
and $T' = T -\{x,y\}$.
Let $G'$ be obtained from $G - T'$ by adding 
a new vertex $z$ and two new edges edge $xz$ and $yz$
{\em (see Figure \ref{fig:2a})}.   
Suppose that   $v(T') = 3k + 2$, where $k \ge 1$ is an integer.
Then
$\lambda (G') \ge v(G')/4 \Rightarrow \lambda (G) \ge v(G)/4$.
\es

\bp 
Clearly $v(G) = v(G') + v(T') - 1$.
Since $v(T') = 3k + 2$ and $k \ge 1$, clearly 
$\lambda (T') = k$.
Let ${\cal P}'$ be a maximum $\Lambda $-packing in $G'$, 
and so $\lambda (G') = |{\cal P}'|$.
Let $xst$ be a subpath of $T$.
\\[1ex]
${\bf (p1)}$ Suppose that $\{xz,yz\} \cap E({\cal P}') = \emptyset $.
Then $\lambda (G) \ge \lambda (G') + \lambda (T') = \lambda (G') + k$.
Now  
$\lambda (G') \ge v(G')/4 \Rightarrow 
\lambda (G) \ge v(G')/4 + k = (v(G) - v(T') +1)/4 + k 
= v(G)/4 + k - (3k+1)/4 = v(G)/4 + (k-1)/4 \ge v(G)/4$.
\\[1ex]
${\bf (p2)}$ Suppose that $|\{xz,yz\} \cap E({\cal P}')| = 1$,
say $rxz$ is a 2-edge path in ${\cal P}'$ for some 
$r \in V(G - T')$.
Let 
${\cal P}''$ be a maximum $\Lambda $-packing in 
$T'' = T' - s$, and
${\cal P} = {\cal P}' \cup {\cal P}''- rxz \cup rxs $. 
Then ${\cal P}$ is a $\Lambda $-packing in $G$, and so 
$\lambda (G) \ge |{\cal P}| = |{\cal P}'| + |{\cal P}''|  = 
\lambda (G') + \lambda (T'')$. 
Since $v(T'') = v(T') - 1 = 3k + 1$, we have $\lambda (T'') = k$. 
Therefore, as in ${\bf (p1)}$,  
$\lambda (G') \ge v(G')/4 \Rightarrow 
\lambda (G) \ge v(G')/4 + k \ge v(G)/4$.
\\[1ex]
${\bf (p3)}$ Now suppose that $xz,yz \in E({\cal P}')$, 
i.e. $xz'y$ is a 2-edge path in ${\cal P}'$.  Let 
${\cal P}''$  be a (unique) maximum $\Lambda $-packing 
in $T'' = T' - \{s,t\}$.
Let ${\cal P} = {\cal P}' \cup {\cal P}''- xzy \cup xst$. 
Then ${\cal P}$ is a $\Lambda $-packing in $G$, and so
$\lambda (G) \ge |{\cal P}| = |{\cal P}'| + |{\cal P}''|  = 
\lambda (G') + \lambda (T'')$. 
Since $v(T'') = v(T') - 2 = 3k$, we have $\lambda (T'') = k$. 
Therefore again, as in ${\bf (p1)}$,  
$\lambda (G') \ge v(G')/4 \Rightarrow 
\lambda (G) \ge v(G')/4 + k  \ge v(G)/4$.
\ep

\bs
\label{Dleaf} 
Let $e \in E(G)$.
Suppose that $G - e = A \cup B$, where 
$A$ and $B$ are disjoint subgraphs of $G$ and $v(A) \ge 3$. 
Suppose also that $\lambda (A) = \lfloor v(A)/3 \rfloor$.
Let $G' = B$.  
Then the implication 
$\lambda (G') \ge v(G')/(s+1) \Rightarrow \lambda (G) \ge 
v(G)/(s+1)$ 
holds for $s = 3$ provided $v(A) \ne 5$, as well as
for every integer $s \ge 4$.
\es

\bp
Clearly $v(G) = v(A) + v(B)$
and $\lambda (G) \ge \lambda (\grave{A}) + \lambda (B) $. Therefore 
$\lambda (B) \ge v(B)/(s+1) \Rightarrow 
\lfloor v(A)/3 \rfloor + \lambda (G) \ge v(B)/(s+1)$.
We claim that if $s \ge 3$, then   
$\lfloor v(A)/3 \rfloor \ge  v(A)/(s+1)$.
This is obviously true for $s \ge 4$ and also true for $s = 3$ because $v(A) \ne 5$. 
Therefore 
$\lambda (G) \ge v(A)/(s+1) + v(B)/(s+1) \ge v(G)/(s+1)$.
\ep

\bs
\label{5-leaf2} 
Let $w = ab \in E(G)$ be a thread of $G$, $a \ne b$, 
$d(a,G) = d(b,G) = 3$, $G - w = A \cup B$, where 
$A$ and $B$ are disjoint subgraphs of $G$, $a \in V(A)$, and
$b \in V(B)$.  
Let $x_iT_ib$, $i = 1,2$, be the two different threads of $G$ distinct from $awb$, 
$bz_i \in E(T_i)$, $\grave{A} = G - (B - b)$ 
(i.e. $\grave {A}  = A \cup awb$ and 
$A' = \grave {A} \cup \{bz_1,bz_2\}$.
Suppose that $v(\grave {A}) = 6$ and 
$\lambda (\grave {A}) = 2$.
\\[0.7ex]  
$(a1)$ If $e(T_i) = 2$ for $i = 1,2$, then let 
$x_iT_ib = x_iz_ib$
{\em (see Figure \ref{fig:3a})}.
\\[0.7ex]  
$(a1.1)$ If $x_1 \ne x_2$, then let $G' = G - A'$.
\\[0.7ex]  
$(a1.2)$ If $x_1 = x_2$, then let $G'$ be obtained from
$G - A'$ by adding a new edge $u = z_1z_2$
{\em (see Figure \ref{fig:4a})}.
\\[0.7ex]  
$(a2)$
If $e(T_1) = 3$, $e(T_2) = 2$, then
let $x_1T_1b = x_1y_1z_1b$, $x_2T_2b = x_2z_2b$.
\\[0.7ex]  
$(a2.1)$ If $x_1 \ne x_2$, then
let $G'$ be obtained from $G - A'$ by adding a new edge 
$u = y_1x_2$ {\em (see Figure \ref{fig:5a})}.
\\[0.7ex]  
$(a2.2)$ If $x_1 = x_2$, then
let $G'$ be obtained from $G - A'$ by adding a new edge 
$u = y_1z_2$ {\em (see Figure \ref{fig:6a})}.
\\[0.7ex]  
$(a3)$ 
If $e(T_1) = 3$, $e(T_2) = 3$, then
let $x_1T_1b = x_1y_1z_1b$, $x_2T_2b = x_2y_2z_2b$ 
and $G'$ be obtained from $G - A'$ by adding a new edge 
$u = y_1y_2$  {\em (see Figure \ref{fig:7a})}.
\\[0.7ex]
\indent
Then $G \in {\cal G}_2^3 \Rightarrow G' \in {\cal G}_2^3$ and
$\lambda (G') \ge v(G')/4 \Rightarrow \lambda (G) \ge v(G)/4$.
\es





\bp It is easy to check that 
$G \in {\cal G}_2^3 \Rightarrow G' \in {\cal G}_2^3$. 
We prove the second claim.
Since  $v(\grave{A}) = 6$ and
$\lambda (\grave{A}) = 2$, clearly 
$v(A') = 8$ and $\lambda (A') = 2$, respectively. 
Let ${\cal P}'$ be a maximum $\Lambda $-packing in $G'$. 
\\[1ex]
${\bf (p1)}$ Suppose that $u \in E(G') \Rightarrow u \not \in E({\cal P}')$. In particular, our assumption holds in case $(a1.1)$. 
  In cases $(a1.2)$ and $(a2.2)$, $u$ belongs to a triangle-leaf of $G'$.   
Hence there exists a maximum $\Lambda $-packing in $G'$ avoiding $u$. 
Therefore we can assume that ${\cal P}'$ avoids $u$, 
and so our assumption holds  in these cases as well.
In cases $(a1.1)$ and (a2.2)  $v(G) = v(G') + v(A')$ and $\lambda (G) \ge \lambda (G') + \lambda (A')$.
Therefore 
$\lambda (G') \ge v(G')/4 \Rightarrow 
\lambda (G) \ge v(G')/4 + \lambda (A') = (v(G) - 8)/4 + 2 = v(G)/4$.
In case $(a1.2)$ $v(G) = v(G') + v(\grave{A})$ and 
$\lambda (G) \ge \lambda (G') + \lambda (\grave{A})$.
Therefore 
$\lambda (G') \ge v(G')/4 \Rightarrow 
\lambda (G) \ge v(G')/4 + \lambda (A') = (v(G) - 6)/4 + 2 > v(G)/4$.
\\[1ex]
${\bf (p2)}$ Suppose that $u \in E(G')$, $u  \in E({\cal P}')$,  
$e(T_1)  = 3$, $e(T_2) = 2$, and $x_1 \ne x_2$, and so 
$(a2.1)$ holds and $u = y_1x_2$.
\\[1ex]
${\bf (p2.1)}$ Suppose that
$x_1y_1\not \in E({\cal P}')$,
say $y_1x_2c_2 \in {\cal P}'$, 
where $c_2 \in V(G' - \{x_2,y_1\})$.
Let ${\cal P}_1$ be a maximum $\Lambda $-packing of 
$A_1 = (A' - z_2) \cup y_1z_1$
and  ${\cal P} = ({\cal P}' \cup {\cal P}_1 - y_1x_2c_2) \cup z_2x_2c_2$.
Then ${\cal P}$ is a $\Lambda $-packing of $G$,
$|{\cal P}| = |{\cal P}'| + |{\cal P}_1|$, and 
$\lambda (G) \ge |{\cal P}| = \lambda (G') +
\lambda (A_1)$. 
Obviously $\lambda (A_1) = \lambda (A') = 2$.
Therefore, as in  ${\bf (p1)}$, 
$\lambda (G') \ge v(G')/4 \Rightarrow v(G)/4$.
\\[1ex]
${\bf (p2.2)}$ Suppose that $x_1y_1\in E({\cal P}')$,
i.e. $x_1y_1x_2 \in {\cal P}'$.
Let ${\cal P}_2$ be a maximum $\Lambda $-packing of 
$A_2 = A' - z_1$ and  
${\cal P} =( {\cal P}' \cup {\cal P}_1 \cup x_1y_1z_1) - x_1y_1x_2$.
Then ${\cal P}$ is a $\Lambda $-packing of $G$,
$|{\cal P}| = |{\cal P}'| + |{\cal P}_2|$, and 
$\lambda (G) \ge |{\cal P}| = \lambda (G') + \lambda (A_2)$. 
Obviously $\lambda (A_2) = 2$.
Therefore, as in  ${\bf (p1)}$, 
$\lambda (G') \ge v(G')/4 \Rightarrow v(G)/4$. 
\\[1ex]
${\bf (p3)}$ Suppose that $u  \in E({\cal P}')$, 
$e(T_1)  = 3$, $e(T_2)  = 3$, and so $(a3)$ holds and
$u = y_1y_2$.
Because of symmetry, we can assume that $x_1y_1y_2 \in {\cal P}'$.
Let ${\cal P}_3$ be a maximum $\Lambda $-packing of 
$A_3 = (A' - z_1) \cup y_2z_2$
and  ${\cal P} = ({\cal P}' \cup {\cal P}_3 - x_1y_1y_2) \cup x_1y_1z_1$.
Then ${\cal P}$ is a $\Lambda $-packing of $G$,
$|{\cal P}| = |{\cal P}'| + |{\cal P}_3|$, and 
$\lambda (G) \ge |{\cal P}| = \lambda (G') + \lambda (A_3)$. 
Obviously $\lambda (A_3) = 2$. 
Therefore, as in  ${\bf (p1)}$, 
$\lambda (G') \ge v(G')/4 \Rightarrow v(G)/4$.  
\ep

\bs
\label{Dstar} 
Let $G$ be a graph, 
$a \in V(G)$, $d(a,G) = d$, 
$\{aT_ix_i: i \in I\}$ be the set of threads of $G$ with 
a common vertex $a$, and 
$S_a = \cup \{T_i - x_i:  i \in I\}$ $($possibly $x_i = x_j$ when $i \ne j$$)$.
Let $s \ge 3$ be an integer and suppose that
$3 \le d \le s$.
Suppose that $e(T_i) \le 3$ for every $i \in I$.
Let $I = I_1 \cup I_2\cup I_3$ be a partition of $I$ such that
$I_1 \ne I$ and 
$(1)$ $i \in I_1$ implies $T_i$ is a triangle,
$(2)$ $i \in I_2$ implies $T_i$ is a 2-edge path, say $x_iz_ia$,  and 
$(3)$ $i \in I_3$ implies  $T_i$ is a 3-edge path, say  $x_iy_iz_ia$.

Let ${\cal T}_k = \{aT_ix_i: i \in I_k\}$ and 
$X_k = \{x \in V(G): aTx \in {\cal T}_k\}$, i.e. $X_k$ is the set of end-vertices of threads in ${\cal T}_k$ that are different from $a$, $k \in \{2,3\}$. Let $t_x$  be the number of threads in ${\cal T}_2$
ending at $x$. Let $X'_2 = \{x \in X_2 \setminus X_3: t_x \ge 2, d(x) = t_x +1\}$. If $x \in X'_2$, then let $L_x$ be a $($unique$)$ thread which ends at $x$ and is not in ${\cal T}$, 
$l_x = e(L_x)$, and $L'_x$ be the subpath of $L_x$ which ends at $x$ and has $l_x$ vertices. 
Suppose that each $l_x \le 3$ and if $|I_3|$ is odd, 
then $I_2 \ne \emptyset $.
Let $G'$ be obtained from $G$ as follows.
\\[0.7ex]
$(a1)$ Suppose that $I_3$ is even. 
Let $S'_a = (S_a - \cup \{y_i: i \in I_3\}) \cup 
(\cup \{L'_x: x \in X'_2\})$.
Partition $I_3$ into pairs $\{i,i'\}$
and put 
$G' = (G \cup \{y_iy_{i'}: i \in I_3\}) - S'_a$
{\em (see Figure \ref{fig:9a})}.
\\[0.7ex]
$(a2)$ Suppose that $I_3$ is odd, and so 
$I_3 \ne \emptyset $.
Choose  $r \in I_3$ and $j \in I_2$  and  partition 
$I_3 - r$ into pairs $\{i,i'\}$.
{\em (Since by our assumption, $I_3$ is odd, $I_2 \ne \emptyset $, and so such $j$ exists.)}
\\[0.7ex]
$(a2.1)$ If $x_r \ne x_j$, then 
let $S'_a = (S_a  - (\cup \{y_i: i \in I_3 - r\} \cup z_j)) 
\cup (\cup \{L'_x: x \in X'_2 - x_j \})$
 and  
$G' = G \cup x_rz_j \cup \{y_iy_{i'}: i \in I_3 - r\} - S'_a$
{\em (see Figure \ref{fig:10a})}. 
\\[0.7ex]
$(a2.2)$ If $x_r = x_j$, then 
let $S'_a = (S_a  - (\cup \{y_i: i \in I_3\} \cup z_j))
(\cup \{L'_x: x \in X'_2 - x_j\})$     
and 
$G' = G \cup y_rz_j \cup \{y_iy_{i'}: i \in I_3 - r\} - S'_a$
{\em (see Figure \ref{fig:11a})}. 
\\[0.7ex]
\indent
$(${\em We say that {\em $G'$ is obtained from $G$ by 
a reduction of the vertex star $S_a$}.}$)$ 
\\[0.7ex]
\indent
Then 
$G \in {\cal G}^s_2 \Rightarrow G' \in {\cal G}^s_2$ and
$\lambda (G') \ge v(G')/(s+1) \Rightarrow \lambda (G) \ge v(G)/(s+1)$.
\es




\bp We can assume that $G$ is a connected graph.
It is easy to check that 
$G \in {\cal G}^s_2 \Rightarrow G' \in {\cal G}^s_2$.
We prove the second claim.
Let ${\cal P}'$ be a maximum $\Lambda $-packing of $G'$, and so $|{\cal P'}| = \lambda (G')$.
Obviously, $v(G) = v(G') + v(S'_a)$. Let $|X'_2| = d'$.
Obviously $d' \le d \le s$.
\\[1ex]
${\bf (p1)}$ Consider case $(a1)$ with $I_3 \ne  \emptyset $, and so $|I_3|$ is even. 
Then $v(G) = v(G') + v(S'_a)$ and  
$v(S'_a) = d + 1 + \sum \{l_x: x \in X'_2\}$.
Let ${\cal P}$ be obtained from ${\cal P}'$ as follows. 
If edge $y_iy_{i'}$ belongs to a
2-edge path $L' \in {\cal P}'$ with the center, say 
$y_k$ in $\{y_i,y_{i'}\}$, then replace
$L'$ in ${\cal P}'$ by the 2-edge path 
$L = (L' - y_{k'}) \cup \{z_k, y_kz_k\}$, where $\{k,k'\} = \{i,i'\}$.
Then $|{\cal P}| = |{\cal P'}| = \lambda (G')$.
Let $S''_a = S'_a - V({\cal P})$.
Then ${\cal P}$ is a $\Lambda $-packing in $G$ and
 $\lambda (G) \ge |{\cal P}| + \lambda (S''_a)$.
 and $\lambda (S''_a) = 1 + d'$.
 Since  $d' \le d \le s$ and each $L_x \le 3$, we have $v(S'_a) \le s+1 + 3d'$. Since $s \ge 3$,clearly 
$\lambda (S'_a) = 1 + d' >  1 + 3d'/(s+1) \ge v(S'_a)/(s+1)$.
Therefore
$\lambda (G') \ge v(G')/(s+1) \Rightarrow 
\lambda (G) \ge \lambda (G') + \lambda (S''_a) \ge
v(G')/(s+1) + v(S'_a)/(s+1)= v(G)/(s+1)$.
\\[1ex]
${\bf (p2)}$ Consider case $(a2.1)$, i.e.  $I_3$ is odd and $x_r \ne x_j$. 
Then $v(G) = v(G') + v(S'_a)$ and  
$v(S'_a) = d + 1 + \sum \{l_x: x \in X'_2\}$.
Let ${\cal P}$ be obtained from ${\cal P}'$ as follows.
If edge $y_iy_{i'}$ in $G'$, where $i \ne r$, belongs to 
a 2-edge path 
$L' \in {\cal P}'$ with the center, say $y_k$ $y_k$ in $\{y_i,y_{i'}\}$, then, as in ${\bf (p1)}$, replace $L'$ in 
${\cal P}'$ by the 2-edge path 
$L = (L' - y_{k'}) \cup \{z_k, y_kz_j\}$, where $\{k,k'\} = \{i,i'\}$. 
If edge $x_rz_j$ in $G'$ belongs to a 2-edge path $P' \in {\cal P}'$ and $P' \ne x_rz_jx_j$ ($P' = x_rz_jx_j$), then
replace $P'$ in ${\cal P}'$ by 2-edge path 
$P = (P' -z_j) \cup \{y_r, x_ry_r\}$
(respectively, by 2-edge path $P = x_ry_rz_r$).
Then $|{\cal P}| = |{\cal P'}| = \lambda (G')$.
Let $S''_a = S'_a - V({\cal P})$.
Then ${\cal P}$ is a $\Lambda $-packing in $G$ and
 $\lambda (G) \ge |{\cal P}| + \lambda (S''_a)$. 
As in  ${\bf (p1)}$,
$\lambda (S''_a) = 1 + d'\ge v(S'_a)/(s+1)$.
Therefore
$\lambda (G') \ge v(G')/(s+1) \Rightarrow 
\lambda (G) \ge \lambda (G') + \lambda (S''_a) 
\ge v(G')/(s+1) + v(S'_a)/(s+1)= v(G)/(s+1)$.
\\[1ex]
${\bf (p3)}$ Consider case $(a2.2)$, i.e. $I_3$ is odd and 
$x_r = x_j$. 
Then $v(G) = v(G') + v(S'_a)$ and  
$v(S'_a) = d + \sum \{l_x: x \in X'_2\}$.
Since $y_rz_j$ belongs to a triangle-bock in $G'$, there exists
a maximum $\Lambda $-packing ${\cal P}'$ in $G'$ that avoids $y_rz_j$.
Let ${\cal P}$ be obtained from ${\cal P}'$ as follows.
If edge $y_iy_{i'}$ in $G'$, where $i \ne r$,  belongs to 
a 2-edge path $L' \in {\cal P}'$ with the center, say  $y_k$ in $\{y_i,y_{i'}\}$, then, as in ${\bf (p1)}$, replace $L'$ in ${\cal P}'$ by the 2-edge path 
$L = (L' - y_{k'}) \cup \{z_k, y_kz_j\}$, where $\{k,k'\} = \{i,i'\}$. 
Then ${\cal P}$ is a $\Lambda $-packing in $G$
and $|{\cal P}| = |{\cal P'}| = \lambda (G')$.
Let $S''_a = S'_a - V({\cal P})$.
Then $\lambda (G) \ge |{\cal P}| + \lambda (S''_a)$.
It is easy to see, that $\lambda (S''_a) = 1 + d'$ and 
$v(S'_a) \le d + 3d' \le s+1 + 3d'$.
Hence, as in  ${\bf (p1)}$, 
$\lambda (S''_a) = 1 + d' \ge v(S'_a)/(s+1)$.
Therefore
$\lambda (G') \ge v(G')/(s+1) \Rightarrow 
\lambda (G) \ge \lambda (G') + \lambda (S''_a) \ge
v(G')/(s+1) + v(S'_a)/(s+1)= v(G)/(s+1)$.
\ep   
\bs
\label{2C5}
Let $a \in V(G)$, $d(a,G) = 3$, 
$aT_ix_i$ be the three threads ending at $a$.
Suppose that $x_j$ belongs to a 5-leaf $X_j$, $e(T_j) = 1$  for $j = 1,2$.
Let $G'$ be obtained from $G - (X_1\cup X_2 - \{x_1,x_2\})$ by adding a new edge $x_1x_2$
{\em (see Figure \ref{fig:8a})}. 
Then 
$\lambda (G') \ge v(G')/4 \Rightarrow \lambda (G) \ge v(G)/4$.
\es

\bp Let $X = X_1\cup X_2 - \{x_1,x_2\}$.
Obviously, $v(G) = v(G') + v(X)$, $v(X) = 8$,
$\lambda (X) = 2$, and 
$\lambda (G) \ge  \lambda (G') + \lambda (X)$.
Therefore 
$\lambda (G') \ge v(G')/4 \Rightarrow 
\lambda (G) \ge  v(G')/4 + \lambda (X) = (v(G) - 8)/4 + 2 = v(G)/4$.
\ep

\section{Constructions of extremal graphs}
\label{constr}
\indent

In this section we give constructions providing infinitely many connected graphs for which  the bounds in 
{\bf \ref{Main1}}, {\bf \ref{DMain1}}, and {\bf \ref{main1}} are attained.
\\[1ex]
\indent
Let $s \ge 1$ be an integer and ${\cal T}_s$ denote the set of trees $T$ such that every non-leaf vertex in $T$ has degree $s$. Obviously $T \in {\cal T}_1 \Rightarrow v(T) = 2$.

Let $k \ge 1$ be an integer and $T_k$ be obtained from $T \in {\cal T}_s$ by subdividing 
each non-dangling edge of $T$ 
with $k$ vertices and each dangling edge with $k-1$ vertices.
Let ${\cal T}_{s,k} = \{T_k: T \in {\cal T}_s\}$.

A $\Lambda _k$-{\em packing} $N$ is a subgraph of $G$ such that every component
of $N$ is a $k$-edge path.
Let $\lambda _k(G)$ denote the maximum number of disjoint $k$-edge paths in $G$.
Obviously $\lambda _1(G) = \tau _1(G)$ is the size of a maximum matching  in $G$,
$\lambda(G) = \lambda _2(G) = \tau _2(G)$, and 
$\lambda _k(G) \le \tau _k(G)$.

The following result provides  infinitely many connected graphs for which the bound in {\bf \ref{main1}} is attained.
\bs
\label{lambda(T)}
Let $k \ge 1$ and $s \ge 1$ be integers and $T \in {\cal T}_s$. 
Then 
$\lambda _k(T_k) = \tau _k(T_k) = (v(T_k) - k)/(sk - k +1)$.
\es

\bp  
We prove our claim on $\lambda _k(T_k)$ by induction on $v(T)$.
The proof on $\tau _k(T_k)$ is similar.
If $s = 1$, then $v(T) = 2$, and so our claim is true. 
So let $s \ge 2$. Then $v(T) \ge s+1$.
If $v(T) = s+1$,  then $T = K_{1,s+1}$, and so our claim is true. 
So let $v(T) > s+1$.
Let $L^x$ denote the set of leaves in $T$ adjacent to 
a vertex $x$.
Since $T \in {\cal T}_s$ and $v(T) > s+1$, 
clearly $T$ has a vertex $z$ with $|L^z| = s-1$.
Let $T' = T - L^z$.
Then  $v(T') = v(T) -(s-1)$ and $T' \in {\cal T}_s$. 
Let $Y^z$ be the subgraph of $T$ induced by $z \cup L^z$.
Clearly $Y^z$ is isomorphic to $K_{1,s-1}$. 
Let $Y^z_k$ be obtained from $Y^z$ by subdividing each edge with $k-1$ vertices. Then $T'_k = T_k - Y^z_k$. 
Obviously
 $\lambda _k(Y^z_k) = 1$, $v(T_k) = v(T'_k) + v(Y^z_k)$,
$v(Y^z_k) = sk - k +1$,
and $T_k$ has a maximum $\Lambda _k$-packing ${\cal P}$ 
that contains exactly one $k$-edge path from $Y^z_k$. 
Hence 
$\lambda _k(T_k) = \lambda _k(T'_k) + \lambda _k(Y^z_k)$.
Since 
\\[0.3ex]
$v(T') < v(T)$, by the induction hypothesis,
$\lambda _k(T'_k)) = (v(T'_k) - k)/(sk - k +1)$. 
Then
\\[0.3ex]
$\lambda _k(T_k)) = (v(T'_k) - k)/(sk - k +1) + 
 \lambda _k(Y^z_k) = (v((T_k)) - (sk - k +1) - k)/(sk - k +1) + 1 =  
\\[0.3ex] 
 (v(T_k) - k)/(sk - k +1)$.
\ep
\\

Let $Y$ denote the graph obtained from three disjoint cycles
$A_i$ of five vertices by adding a new vertex $a$ and three edges
$\{aa_1,aa_2,aa_3\}$, where $a_i \in V(A_i)$, $i = 1,2,3$.
It is easy to see that $\lambda (Y) = v(Y)/4 = 4$.
Here is a more general construction of extremal graphs
which shows, in particular, that there are infinitely many connected graphs for which the bounds in {\bf \ref{Main1}} and in {\bf \ref{DMain1}} are attained.  
\bs
\label{2subd}
Let $H$ be a graph with possible loops and parallel edges 
such that each vertex of $H$ is of degree at least three and 
at most $s$. 
Let $H_k$ be a graph obtained from $H$ by subdividing 
every  edge with exactly  $k \ge 1$ vertices. 
Then $\lambda _k(H_k) = \tau _k(H_k) = v(H)$.
Moreover, $\lambda (H_2) =  v(H) \ge v(H_2)/(s+1)$
and the equality holds if and only if $H$ is an $s$-regular graph.
\es

\bp Obviously $\lambda _k(H_k) \le \tau _k(H_k)$.
Let ${\cal P}$ be a $\Lambda _k$-packing in $H_k$.  
Then each $k$-edge tree in ${\cal P}$ contains a vertex 
from $V(H)$.
Therefore $\lambda _k(H) \le \tau _k(H_k) \le v(H)$.
We can assume that $H$ is connected. 
Let $T$ be a spanning tree of $H$.
Obviously there is an edge $e = xy \in E(H) - E(T)$.
Let $D = T \cup e$ and 
$D_x$ be the directed graph obtained from $D$ by directing each edge of 
$T \subset D$ towards $x$ and by directing $e$ from $x$ to $y$.
For $z \in D_x$, let $p(z)$ denote the edge in $D_x$ with the tail $z$.
Clearly $p: V(D_x) \to E(D_x)$ is a bijection.
Let  $\ddot{D}$ be obtained from $D$ by subdividing each edge of $D$ with exactly $k$ vertices. Then every edge $e$ in $D$ is replaced in $\ddot{D}$ by a thread which we denote by
$L(e)$. Let $p_k(z)$ denote the $k$-edge path in 
thread $L(p(z)$ containing $z$.
We can assume that $\ddot{D} \subseteq H_k$.
Let ${\cal P} = \{p_k(z): z \in V(D_x)\}$.
Then clearly ${\cal P}$ is  a $\Lambda _k$-packing in 
$\ddot{D}$, and therefore in $H_k$, and 
$|{\cal P}| = v(T) = v(H)$. 
Therefore $\lambda _k(H) = \tau _k(H_k) = |{\cal P}| = v(H)$.

Obviously, $2e(H) \le sv(H)$ and $v(H_2) =  v(H) + 2e(H) \le (s+1)v(H)$.
Therefore $\lambda (H_2)  = v(H) \ge  v(H_2)/(s+1)$ and 
the equality holds if and only if $H$ is an $s$-regular graph.
\ep

\section{Proof of Theorem \ref{Main1}}
\label{Proof}

\indent

First we will establish some properties of 
a hypothetical minimum counterexample to theorem {\bf \ref{Main1}}. 
We will use these properties 
to proof theorem {\bf \ref{Main1}} by showing that no counterexamples exist.
\\

In all claims below, except for {\bf \ref{Main1}}, we assume that 
$G$ is a connected graph satisfying the following conditions:
\\[1ex]
$(a1)$ $G$ is an ${\cal F}^3_2$-minimal graph,
\\[1ex]
$(a2)$ $\lambda (G) < v(G)/4$, and
\\[1ex]
$(a3)$ $G$ has the minimum number of vertices among all graphs satisfying  $(a1)$ and $(a2)$.

\bs
\label{nothread>3}
Let $aTb$ be a thread of $G$.
If neither $a$ nor $b$ belongs to a 5-leaf, then 
$e (T) \in \{2,3\}$.
\es

\bp (uses {\bf \ref{Dthread1}} and {\bf \ref{thread2}}).
Suppose, on the contrary, that $e(T) \not \in \{2,3\}$.
Obviously $e(T) \ge 1$.
If $T$ has exactly one edge $u$, then $G - u \in {\cal F}^3_2$, 
and so $G$ is not  ${\cal F}^3_2$-minimal, a contradiction.
Therefore we have exactly one of the following two cases:
\\[1ex]
$(c1)$  $e(T) \in \{3k+1, 3k + 2\}$, where $k \ge 1$ and
\\[0.7ex]
$(c2)$ $e(T) = 3k$, where $k \ge 2$.
\\[1ex]
${\bf (p1)}$
Suppose that $e(T) \in \{3k+1, 3k + 2\}$, where $k \ge 1$.
Let, as in {\bf \ref{Dthread1}}, $G' = G - (T -\{a,b\})$.
Since $G \in {\cal F}^3_2$ and neither $a$ nor $b$ belongs to a 5-leaf, we have $G' \in {\cal F}^3_2$.
Clearly $v(G') < v(G)$.
By the minimality of $G$, $\lambda (G') \ge v(G')/4$.
Then by {\bf \ref{Dthread1}}, $\lambda (G) \ge v(G)/4$, a contradiction.
\\[1ex]
${\bf (p2)}$ Now suppose that $e(T) = 3k$, where $k \ge 2$.
Let $G'$ be the graph defined in {\bf \ref{thread2}}.
Since $G \in {\cal F}^3_2$, obviously $G' \in {\cal F}^3_2$.
Clearly $v(G') < v(G)$.
By the minimality of $G$, $\lambda (G') \ge v(G')/4$.
Then by {\bf \ref{thread2}}, $\lambda (G) \ge v(G)/4$, a contradiction.
\ep

\bs
\label{no5-leaf}
$G$ has no 5-leaves.
\es

\bp (uses 
{\bf \ref{Dleaf}}, 
{\bf \ref{5-leaf2}}, 
{\bf \ref{2C5}}, and
{\bf \ref{nothread>3}}).
Suppose, on the contrary, that $G$ has a 5-leaf $X_1$.
Let $aT_1x_1$ be the stem of $X_1$, i.e. a (unique) thread 
in $G$ such that $x_1 \in V(X_1)$.
If $a$ belongs to a 5-leaf or a cycle-leaf $A$, then 
$G = A\cup T_1 \cup X_1$ and $G$ has a Hamiltonian path.
Hence $\lambda (G) = \lfloor v(G)/3 \rfloor$. 
Obviously $v(G) > 5$. Therefore
$\lambda (G) = \lfloor v(G)/3 \rfloor \ge v(G)/4$, 
a contradiction.
So we assume that $a$ belongs to neither a 5-leaf nor a cycle-leaf.
\\[1ex]
${\bf (p1)}$
Suppose that $e(T_1) \ge 2$. 
Let $\grave{X}_1 = X_1 \cup (T_1 - a)$ and, 
as in {\bf \ref{Dleaf}}, $G' = G - \grave{X}_1$.
Then $\lambda (\grave{X}_1) = \lfloor v(\grave{X}_1/3\rfloor $.
Since $G \in {\cal F}^3_2$ and $a$ belongs to no 5-leaf,
we have $G' \in {\cal F}^3_2$.
Clearly $v(G') < v(G)$.
By the minimality of $G$, $\lambda (G') \ge v(G')/4$.
Then by {\bf \ref{Dleaf}}, $\lambda (G) \ge v(G)/4$, 
a contradiction.
\\[1ex]
${\bf (p2)}$
Now suppose that $e(T_1) = 1$.
Let $aT_2x_2$ and $aT_3x_3$ be the threads containing $a$ and distinct from $aT_1x_1$.
Since $a$ belongs to no cycle-leaf, clearly $T_2 \ne T_3$.
If $x_i$ belongs to no 5-leaf, then by {\bf \ref{nothread>3}}, 
$e(T_i) \in \{2,3\}$.
\\[1ex]
${\bf (p2.1)}$
Suppose that neither $x_2$ no $x_3$ belongs to a 5-leaf.
Then $e(T_i) \in \{2,3\}$ for $i = 1,2$.
Let $G'$ be the graph defined in {\bf \ref{5-leaf2}}, 
where $\grave{A} = X_1 \cup T_1$. 
Then 
$v(\grave{A}) = 6$ and $\lambda ((\grave{A}) = 2$.
By  {\bf \ref{5-leaf2}}, $G \in {\cal F}^3_2$ implies 
$G' \in {\cal F}^3_2$. Clearly $v(G') < v(G)$.
By the minimality of $G$, $\lambda (G') \ge v(G')/4$.
Then by {\bf \ref{5-leaf2}}, $\lambda (G) \ge v(G)/4$, a contradiction.
\\[1ex]
${\bf (p2.2)}$
Now suppose that $x_i$ belongs to a 5-leaf $X_i$ for some 
$i \in \{2,3\}$, say for $i = 2$.
Then by the above arguments on $e(T_1)$, we can assume that   $e(T_2) = 1$.
Let $G'$ be the graph defined in {\bf \ref{2C5}}.
Then $G \in {\cal F}^3_2$ implies $G' \in {\cal F}^3_2$.
Clearly $v(G') < v(G)$.
By the minimality of $G$, $\lambda (G') \ge v(G')/4$.
Then by {\bf \ref{2C5}}, $\lambda (G) \ge v(G)/4$, a contradiction.
\ep

\bs
\label{allthreads=3}
Every path-thread in $G$ has exactly three edges.
\es

\bp (uses 
{\bf \ref{Dstar}}, 
{\bf \ref{nothread>3}}, and 
{\bf \ref{no5-leaf}}).  
By {\bf \ref{no5-leaf}},  $G$ has no 5-leaves.
Let $aT_1x_1$ be a path-thread in $G$.
By {\bf \ref{nothread>3}}, $e(T_1) \in \{2,3\}$.
Suppose, on the contrary, that $e (T_1) = 2$.
If each end-vertex of $T$ belongs to a cycle-leaf, then
obviously $\lambda (G) \ge v(G)/4$. 
Therefore we can assume that one end-vertex,
say $a$, of $T_1$ belongs to no cycle-leaf.

Let $aT_1x_1$, $aT_2x_2$, and $aT_3x_3$ be the three different path-threads ending at $a$.
Since $G$ has no 5-leaves, by {\bf \ref{nothread>3}}, 
$e(T_2) \in \{2,3\}$ and $e(T_3) \in \{2,3\}$. 
Let $x_ix'_i, aa_i \in E(T_i)$ for $i \in \{1,2,3\}$, and so 
$a_i = x'_i$ if and only if $e(T_i) = 2$. 

If  $x_1 = x_2 = x_3 = x$, then $G$ consists of three threads with common end-vertices $a$ and $x$. Since $G \in {\cal F}_2^3$, we have $v(G) \ge 6$. 
Therefore $v(G) \in \{6, 7\}$. 
Obviously $\lambda (G) = 2  = \lceil v(G)/4 \rceil $, a contradiction.
Therefore we assume that $\{x_1, x_2, x_3\}$ has at least two different vertices. Thus we have (up to symmetry) the following three cases:
$(c1)$ $x_1 \ne x_2 = x_3 = x$,
$(c2)$ $x = x_1 = x_2 \ne x_3$, and 
$(c3)$ all $x_i$'s are different.
\\[0.5ex]
\indent
Let $G'$ be a graph obtained from $G$ by a reduction of the vertex star $S_a$ (see {\bf \ref{Dstar}}).
Namely (we remind that in this particular case), $G'$ is obtained from $G$ as follows.  

Consider case $(c1)$, i.e. $x_1 \ne x_2 = x_3 = x$. 
If $\{e(T_2), e((T_3)\} = \{2,3\}$, say 
$\{e(T_2)= 2$ and $e(T_3) = 3$,
then $G' = (G - a_1aa_3) \cup x'_2x'_3$.
If $e(T_2) = e(T_3) = 2$, then 
$G' = G - (T_1 \cup T_2 \cup T_3 \cup Z - x_1)$, where
$Z$ is the set of inner vertices of the thread containing $x$ and distinct from $T_i$, $i \in \{2,3\}$.
If $\{e(T_2) = e((T_3)\} = 3$, then 
$G' = (G - \{a,a_1,a_2,a_3\}) \cup x'_2x'_3$.

Consider case $(c2)$, i.e. $x = x_1 = x_2 \ne x_3$.
If $e(T_2)= 2$ and $e(T_3) = 3$, then we can put
$G' = (G \cup a'_1x_3) - a'_2aa'_3x'_3$.
If
$e(T_2)= 3$ and $e(T_2) = 3$, then we can put
$G' = (G \cup x'_1x'_2) - a'_2aa'_3$.
If $e(T_2) = e(T_3) = 2$,  then 
$G' = G - (T_1 \cup T_2 \cup T_3 \cup Z - x_1)$, where
$Z$ is the set of inner vertices of the thread containing $x$ and distinct from $T_i$, $i \in \{1,2\}$.
If $\{e(T_2) = e((T_3)\} = 3$, then 
$G' = (G \cup x'_2x'_3) - (T_1 \cup T_3 - \{x_1,x_3\})$.

Consider case $(c3)$, i.e.
all $x_i$'s are different.
If $e(T_2) = e(T_3) = 2$, then $G' = G - \{a, a_1,a_2, a_3\}$.
If $\{e(T_2), e(T_3)\} = \{2,3\}$, say $e(T_2) = 2$ and $e(T_3) = 3$, 
then we can put $G' = (G - \{a,a_1,a_3, x'_3\} \cup x'_2x_3$.
If $e(T_2) = e(T_3) = 3$, then 
$G' = (G - \{a,a_1,a_2,a_3\}) \cup x'_2x'_3$.
\\[0.5ex]
\indent
By  {\bf \ref{Dstar}}, 
$G \in {\cal F}^3_2 \subset {\cal G}^3_2 \Rightarrow G' \in {\cal G}^3_2$.
Since {\bf \ref{no5-leaf}},  $G$ has no 5-leaves, we have 
$G' \in {\cal F}^3_2$.
Clearly $v(G') < v(G)$.
By the minimality of $G$, $\lambda (G') \ge v(G')/4$.
Then by {\bf \ref{Dstar}}, $\lambda (G) \ge v(G)/4$, 
a contradiction.
\ep

\bs
\label{3-leaf}
Every cycle-leaf of $G$ is a triangle. 
\es

\bp (uses {\bf \ref{Dleaf}}, {\bf \ref{no5-leaf}}, and 
{\bf \ref{allthreads=3}}). 
Let $A$ be a cycle-leaf of $G$ and $aTb$ be the stem of $A$ with $a \in V(A)$.

Suppose that $v(\grave{A}) \ne 5$. 
Let, as in {\bf \ref{Dleaf}}, $G' = G - \grave{A}$.
Obviously 
$\lambda (\grave{A}) = \lfloor v(\grave{A})/3 \rfloor $.
By {\bf \ref{no5-leaf}}, $G$ has no 5-leaves.
Therefore $G \in {\cal F}^3_2$ implies $G' \in {\cal F}^3_2$.
Clearly $v(G') < v(G)$.
By the minimality of $G$, $\lambda (G') \ge v(G')/4$.
Then by {\bf \ref{Dleaf}}, $\lambda (G) \ge v(G)/4$, 
a contradiction.

Now suppose that $v(\grave{A}) = 5$, i.e. $v(A) + e(T) = 6$.
By {\bf \ref{allthreads=3}}, $e(T) = 3$. Therefore $v(A) = 3$.
\ep
\\

Now we are ready to prove 
\\[1ex]
{\bf \ref{Main1}}~~
{\em Let $G \in {\cal F}^3_2$. 
Then $\lambda (G) \ge v(G)/4$.}
\\[1ex]
\bp (uses 
{\bf \ref{min}},
{\bf \ref{2subd}}, 
{\bf \ref{allthreads=3}}, and 
{\bf \ref{3-leaf}}).
Suppose, on the contrary, that our claim is not true.
Obviously, a triangle is the minimum graph in ${\cal F}^3_2$
and our claim is true for a triangle.
Let $G$ be a vertex minimum counterexample.
By {\bf \ref{min}}, we can assume that $G$ is 
${\cal F}^3_2$-minimal.
In other words, $G$ is a graph such that
$(a1)$ $G$ be an ${\cal F}^3_2$-minimal graph,
$(a2)$ $\lambda (G) < v(G)/4$, and 
$(a3)$ $G$ has the minimum number of vertices among all graphs satisfying 
$(a1)$ and $(a2)$.
By {\bf \ref{allthreads=3}}, each path-thread of $G$ has exactly three edges.
By {\bf \ref{3-leaf}}, each leaf of $G$ is a triangle.
Then by {\bf \ref{2subd}}, $\lambda (G) = v(G)/4$, a contradiction.
\ep
\\

The construction in {\bf \ref{2subd}} provides infinitely many
2-connected graphs (moreover, subdivisions of 3-connected
graphs) $G$ in ${\cal F}^3_2$ such that 
$\lambda (G) = v(G)/4$.

\section{Proof of Theorem \ref{DMain1}} 

\indent

In all claims below, except for {\bf \ref{DMain1}}, we assume that 
$G$ is a connected graph satisfying the following conditions:
\\[1ex]
$(a1)$ $G$ is an ${\cal G}^s_2$-minimal graph, $s \ge 4$,
\\[0.7ex]
$(a2)$ $\lambda (G) < v(G)/(s+1)$, and 
\\[0.7ex]
$(a3)$ $G$ has the minimum number of vertices among all graphs satisfying  $(a1)$ and $(a2)$.

\bs
\label{Dnothread>3}
Let $xTy$ be a thread of $G$, where possibly $x = y$.
Then $e(T) \in \{2,3\}$.
\es

\bp (uses {\bf \ref{Dthread1}}).
Suppose, on the contrary, that $e(T) \not \in \{2,3\}$.
Obviously $e(T) \ge 1$.
If $T$ has exactly one edge $u$, then $G - u \in {\cal G}^s_2$, and so $G$ is not 
${\cal G}^s_2$-minimal, a contradiction.
Therefore $e(T) \ge 4$.
Let, as in {\bf \ref{Dthread1}}, $G' = G - (T - \{x,y\})$.
Obviously $G \in {\cal G}^s_2$ implies $G' \in {\cal G}^s_2$.
Clearly $v(G') < v(G)$.
Therefore by the minimality of $G$, 
$\lambda (G') \ge v(G') /(s+1)$. 
Then by {\bf \ref{Dthread1}}, $\lambda (G) \ge v(G) /(s+1)$,
a contradiction.
\ep

\bs
\label{Dallthreads=3}
Every thread of $G$ has exactly three edges.
\es

\bp (uses {\bf \ref{Dstar}} and {\bf \ref{Dnothread>3}}).
If $T$ is a cycle-thread of $G$, then by {\bf \ref{Dnothread>3}},
$T$ is a triangle, and so $e(T) = 3$.
So let us consider a path-thread $a_1T_1x_1$.
By {\bf \ref{Dnothread>3}}, $e(T_1)\in \{2,3\}$.

Suppose, on the contrary, that $e(T_1) = 2$. 
Let $\{T_i:i \in \{1, \ldots, k\}\}$ be the set of threads of $G$ with
a common end-vertex $a$. 
Let $G'$ be a graph defined in {\bf \ref{Dstar}}.
Then by $G \in {\cal G}^s_2$ implies $G' \in {\cal G}^s_2$.
Clearly $v(G') < v(G)$.
Therefore by {\bf \ref{Dstar}}, the minimality of $G$, 
$\lambda (G) \ge v(G) /(s+1)$. 
Then by {\bf \ref{Dstar}}, $\lambda (G) \ge v(G) /(s+1)$,
a contradiction. 
\ep
\\

Now we are ready to prove  
\\[1ex]
{\bf \ref{DMain1}}~~ 
{\em Let $G \in {\cal G}^s_2$ and $s \ge 4$.
Then $\lambda (G) \ge v(G)/(s+1)$.}
\\[1ex]
\bp (uses 
{\bf \ref{min}}, 
{\bf \ref{2subd}}, and 
{\bf \ref{Dallthreads=3}}).
Suppose, on the contrary, that our claim is not true.
Obviously, a triangle is the minimum graph in ${\cal G}^s_2$
and our claim is true for a triangle.
Let $G$ be a vertex minimum counterexample.
By {\bf \ref{min}}, we can assume that $G$ is 
${\cal G}^s_2$-minimal.
In other words, $G$ is a graph such that
$(a1)$ $G$ be an ${\cal G}^s_2$-minimal graph,
$(a2)$ $\lambda (G) < v(G)/(s+1)$, and 
$(a3)$ $G$ has the minimum number of vertices among all graphs satisfying 
$(a1)$ and $(a2)$.
By {\bf \ref{Dallthreads=3}}, each thread of $G$ has exactly three edges.
Then by {\bf \ref{2subd}}, $\lambda (G) \ge v(G)/(s+1)$, 
a contradiction.
\ep
\\

The construction in {\bf \ref{2subd}} provides infinitely many
2-connected graphs (moreover, subdivisions of 3-connected
graphs) $G$ in ${\cal G}^s_2$ such that 
$\lambda (G) = v(G)/(s+1)$.

\section{Proof of Theorem \ref{main1}}

\indent 

A subtree $B$ of a tree is called a {\em branch} of $T$
if either $B = T$ or $T - B$ is also a tree. If
$B \ne T$, then let $r(B)$ be the vertex of $B$ adjacent to a vertex in $T - B$.
We call $r(B)$ the {\em root} of the branch $B$.

Let, as above,  $\tau _k(G)$ denote the maximum number  of disjoint $k$-edge trees in $G$.
\bs
\label{branch}
Let $s \ge 3$ and $k \ge 1$ be integers, $T$ a tree, 
$T \in {\cal G}^s_1$, and 
$v(T) \ge k+1$.  Then $T$ has a branch $B$ such that 
$v(B) \le (s-1)k + 1$ and $\tau _k(B) = 1$.
\es

\bp
We prove our claim by induction on $v(T)$.
If $v(T) = k+1$, then our claim is obviously true. 
Therefore let $v(T) > k+1$.
A branch $B$ of a tree $T \in {\cal G}^s_1$ is called 
$k$-{\em good} if $v(B) \le (s-1)k + 1$ and $\tau _k(B) = 1$. 
Let $x$ be a leaf of $T$ and $T' = T - x$.
Clearly $T'$ is a tree, $T'\in {\cal G}^s_1$, and  $v(T') < v(T)$.
Therefore by the induction hypothesis,
$T'$ has a $k$-good branch $B'$. 
Let $r' = r(B')$ and $xy \in E(T)$.
Let us assume that $B'$ is a $k$-good branch in $T'$ having the minimum number of vertices. 
If $y \in V(T - (B' - r')$, then $B'$ is a required branch of $T$.
Therefore let $y \in V(B' - r')$. 
Since $\tau _k(B') = 1$ and $B'$ is a vertex minimum 
$k$-good branch in $T'$,
clearly $\tau _k(C') = 0$ for every component $C'$ of $B' - r'$,  
and so $v(C') \le k$. 
Therefore every $k+1$-vertex subtree $S$ in $B'$ 
contains $r'$.
Let $X'$ the component of $B' - r'$ containing $y$,
$X = X' \cup xy$, and $B = B' \cup xy$. 
Clearly $B$ is a branch of $T$ and $r' = r$ is the root of $B$.

Suppose that $v(X') < k$. Then $v(X) = v(X) +1 \le k$, 
$v(B') < (s-1)k + 1$, and so $v(B) = v(B') + 1 \le (s-1)k + 1$.
Since $\tau _k(B') = 1$, clearly $\tau _k(B) \ge 1$, and so $B$ has a $k$-edge subtree. 
Since $v(C) \le k$ for every component $C$ of $B - r$, every subtree of $k$ edges in $B$ contains $r$, and so 
$\tau _k(B) = 1$. Therefore $B$ is a $k$-good branch in $T$.

Now suppose that $v(X') = k$. 
Then $X$ is a $k$-edge subtree of $T$,
$X$ is a branch of $T$, $v(X) \le v(B') (s-1)k + 1$, and 
$\tau _k(X) = 1$. 
Therefore $X$ is a $k$-good branch in $T$. 
\ep
\\

Now we are ready to prove 
\\[1ex]
{\bf \ref{main1}}
{\em Let $s \ge 3$ and $k \ge 1$ be integers,
$G \in {\cal G}^s_1$, and $G$ has no $k$-vertex 
components. 
Then $\tau _k(G) \ge  (v(G) - k)/(sk - k +1)$.}
\\[1ex]
\bp (uses {\bf \ref{branch}}). 
Let $G \in {\cal G}^s_1$. 
We can assume that $G$ is connected and $v(G) \ge k$. 
Since $G$ has no $k$-vertex component, 
$v(G) \ge k+1$. Let $T$ be a spanning tree of $G$. 
Clearly $T \in {\cal G}^s_1$ and 
$\lambda _k(G) \ge \lambda _k(T)$.
Therefore it is sufficient to proof our claim for every tree $T$ in ${\cal G}^s_1$.

We prove our claim for trees by induction on $v(T)$.
If $v(T) = k+1$, then our claim is obviously true. 
So let $v(T) > k+1$. 
By {\bf \ref{branch}}, $T$ has a branch $B$ such that 
$v(B) \le (s-1)k + 1$ and $\tau _k(B) = 1$.
Let $T' = T - B$.
Clearly $T' \in {\cal G}^s_1$, 
$v(T')  
< v(T)$, and
$\tau _k(T ) \ge \tau _k(T') + \tau _k(B)$.
By the induction hypothesis,
$\tau _k(G') \ge  (v(G') - k)/(sk - k +1)$.
Therefore 
$\tau _k(G) \ge (v(G') - k)/(sk - k +1) + 1 =
(v(T) - v(B) - k)/(sk - k +1) + 1$.
Since $v(B) \le (s-1)k + 1$, we have
$\tau _k(G) \ge 
(v(T) - (sk - k + 1) - k)/(sk - k +1) + 1 = (v(G) - k)/(sk - k + 1)$.
\ep
\\

The construction in  {\bf \ref{lambda(T)}}
provides infinitely many trees $T$ in ${\cal G}^s_1$ such that 
$\tau _k(G) = (v(G) - k)/(sk - k +1)$.



\begin{thebibliography}{9}
\bibitem{Khrd} A. Kelmans, Packing $P_k$ in a cubic graph is $NP$--hard 
for $k \ge 3$, manuscript.
\end{thebibliography}
\end{document}